\documentclass[10pt]{amsart}
\theoremstyle{plain}
\usepackage{amsmath, amssymb,graphicx,amscd,color}
\newtheorem{lem}{Lemma}[section]
\newtheorem{prop}[lem]{Proposition}
\newtheorem{cor}[lem]{Corollary}
\newtheorem{dfn}[lem] {Definition}

\newtheorem*{tm}{Theorem}

\newcommand{\co}{{\colon}}
\newcommand{\s}{{\mathfrak{s}}}
\newcommand{\so}{{\overline{\mathfrak{s}}_o}}
\newcommand{\sss}{{\overline{\mathfrak{s}}}}

\newcommand{\x}{{\bf v}}
\newcommand{\ckh}{{CKh}}
\newcommand{\kh}{{Kh}}

\newcommand{\R}{{\bf{R}}}

\newcommand{\Z}{{\bf{Z}}}

\newcommand{\A}{{{\mathcal A }}}

\newcommand{\ts}{{{\thinspace}}}

\title{Khovanov's invariant for closed surfaces}
\author{Jacob Rasmussen}
\address{Princeton University Dept. of Mathematics, Princeton, NJ 08544}
\email{jrasmus@math.princeton.edu}
\thanks{The author was partially supported by an NSF Postdoctoral
  fellowship.}

\begin{document}

\begin{abstract}
We show that the Khovanov-Jacobsson number of an embedded torus in \(\R^4\) is
always \( \pm2\).  
\end{abstract}

\maketitle

\section{Introduction}

In his original paper describing his Jones polynomial homology for
knots \cite{Khovanov}
 --- now known as the Khovanov homology --- Khovanov described
how his construction could be used to define an invariant of
embedded surfaces in \(\R^4\). More precisely, let \( S \subset
\R^3 \times [0,1] \) be a smoothly embedded, orientable cobordism
between knots \(K_0\) and \(K_1\). Given a decomposition of 
\( S \) into  elementary cobordisms, Khovanov defined a
graded map
\( \phi_{S} \colon Kh(K_0) \to Kh(K_1)\) and conjectured that up
to sign, \( \phi_{S}\) was an invariant of \(S
\). This conjecture was subsequently proved by Jacobsson
\cite{Jacobsson} and Khovanov \cite{Khovanov3}, and later 
in a more general form by Bar-Natan \cite{DBN2}.

The Khovanov homology of the empty link is \(\Z\), so this
construction associates a homomorphism
\( \phi_\Sigma : \Z \to \Z \) (well defined up to sign) to a closed,
orientable, smoothly embedded surface \( \Sigma \subset \R^4\).
Defining \(n_\Sigma = |\phi_\Sigma (1)|\), we obtain an invariant of
 \(\Sigma\), which is generally referred to as the Khovanov-Jacobsson
 number.  Since \( \phi_\Sigma\) is a graded map
of degree \(\chi(\Sigma)\), \(n_\Sigma = 0 \) unless \(\chi (\Sigma) =
0\). Thus, the simplest case in which one might have a nonzero
invariant is when \(\Sigma\) is a torus. In \cite{CSM}, 
Carter, Saito, and Satoh showed that \(n_\Sigma = 2
  \) for a large class of embedded tori, namely those which
  can be unknotted by a double point move. In fact, this relation
  holds for all embedded tori in \(\R^4\):

\begin{tm}
If \( \Sigma \subset \R^4\) is a smoothly embedded torus,
 then \( n_\Sigma = 2 \). 
\end{tm}
\noindent The same result has also been obtained by  Tanaka
\cite{Tanaka} using a slightly different method. 

The proof of the theorem is a straightforward extension of the
techniques of \cite{khg}. We briefly recall the setup of that
paper, and refer the reader to it for more details. Let \(D\) be a
planar diagram of a link \(L\). Associated to \(D\), there is a graded
chain complex \( CKh(D) \) whose homology is \(Kh(L)\). In
\cite{ESL2}, Lee defined a related {\it filtered} complex
\(CKh'(D)\). The filtration gives rise to a spectral sequence whose
\(E_1\) term is the complex \(CKh(D)\). This spectral sequence
converges to the homology group \(Kh'(L)\), which she explicitly
calculated. 

To be precise, if \(L\) has \(l\) components, 
then \(Kh'(L) \) has rank \(2^l\) over any field
with characteristic not equal to \(2\). Moreover, Lee gave an explicit
correspondence between orientations on \(L\) and generators of
\(Kh'(L)\). In \cite{khg}, it was observed that these generators are
``projectively canonical,'' in the sense that they are preserved up
to scalar multiplication by the isomorphisms induced by Reidemeister
moves. As we explain in section~\ref{Sec:Gens}, 
it turns out that after a slight modification, 
Lee's generators are actually canonical up to sign.
 This enables us to give a more precise
calculation of the maps induced by cobordisms in Lee's
theory. We then explain how this calculation can
be  used to prove the theorem above. 

The author would like to thank Ivan Smith for helpful conversations
in which this question was raised, and Marc Lackenby and the
organizers of the 2004 Princeton-Oxford topology conference, where
this work was done. 

\section{Cobordisms}
\label{Sec:Cobordisms}

We begin with some generalities on cobordisms, orientations, and
induced maps. Let \(S \subset \R^3 \times[0,1]\) be  a smoothly,
properly embedded, orientable surface with boundary \(L_1 \times {0}
\cup L_2 \times{1} \). We say that \(S\) is a {\it cobordism}  from
\(L_1\) to \(L_2\), and write \(S : L_1 \to L_2\).
Next, suppose that  \(o_1\) and \(o_2\)  are 
 orientations on \(L_1\) and \(L_2\). We say that
  \(o_1\) and \(o_2\) are {\it compatible} if there is an orientation
 \(o\) on \(S\) which induces \(o_1\) on \(L_1\) and 
 the reverse of \(o_2\) on \(L_2\).
(Thus if \(S:L_1 \to L_1\) is the identity
cobordism, \(o_1\) is compatible with itself.) 

\subsection{Induced Maps}
If \(S_1: L_1 \to L_2 \) and \(S_2: L_2 \to L_3\) are cobordisms, 
 we can compose them to get a cobordism \(S:L_1 \to
L_3\). Conversely, it is well known \cite{CarterSaito} that any
cobordism is isotopic to a composition of certain {\it elementary
  cobordisms} which may be represented by fixed local moves on 
planar diagrams. The necessary local moves are given by the
Reidemeister moves, their inverses, and the Morse moves (addition of a
\(0,1\) or \(2\)-handle.)

Following Khovanov \cite{Khovanov}, we want to assign to a
cobordism \(S \co L_1 \to L_2\) an induced map \( \phi_S \co \kh(L_1)
\to \kh(L_2) \) and \(\phi_S'\co \kh'(L_1)
\to \kh'(L_2) \). (Note that this notation differs from that of
\cite{khg}, where the induced map on \(\kh'\) was denoted by
\(\phi_S\).) Since we would like this assignment to be functorial, it
suffices to define the induced maps associated to the various 
elementary cobordisms.

First, suppose \(S: L \to L\) is a Reidemeister move relating two diagrams
\(D\) and \(\tilde{D}\) of \(L\). In \cite{Khovanov}, 
Khovanov 
defined  maps \( \rho_i \co \ckh (D) \to \ckh(\tilde{D}) \) and 
\( \rho_i' \co \ckh (D) \to \ckh(\tilde{D}) \) which induce
isomorphisms on homology. In this case, we define \( \phi_S \) to be \(
\rho_{i*}\). Likewise, if \(S: L \to L\) is the inverse of a
Reidemeister move, \( \phi_S = \rho_{i*}^{-1} \). 

Now suppose \(S
\co D_1 \to D_2 \) is an elementary cobordism associated to a Morse
move. To define \( \phi_S\) in this case, we recall that the chain
complex  \(\ckh(D_1)\) is defined using a certain \(1+1\)-dimensional
TQFT \(\A\). \(\ckh(D_1) \) is generated by elements 
 of the form \(\x \in \A (D_{1,v})\), 
where \(D_{1,v}\) is a complete resolution of the diagram \(D_1\).
\(D_{1,v}\) naturally determines a complete resolution \(D_{2,v}\) of \(D_2\),
and the
cobordism \(S\) induces a cobordism \(S_v \co D_{1,v} \to D_{2,v} \).
We define a map \( \phi \co \ckh(D_1) \to \ckh(D_2) \)
 by setting \( \phi (\x) = \A(S_v)(\x)\) for \(\x \in D_{1,v}\). 
 \( \phi_S\) is defined to be the induced map on homology. With these
 definitions, Jacobsson \cite{Jacobsson} and Khovanov \cite{Khovanov3}
 showed that \( \phi_S\) is well defined up to sign, regardless of how
 \(S\) was  decomposed into elementary cobordisms.

Exactly the same procedure can be used to define an induced map
\(\phi_S'\co \kh'(L_1) \to \kh'(L_2) \) in Lee's theory. Indeed,
Lee's theory is formally identical to Khovanov's, but with the TQFT
\(\A\) replaced by a new TQFT \( \A'\). (See \cite{DBN2} for a
realization of both theories as specializations of a more general,
geometric theory, and for a universal proof that maps induced by
cobordisms are well defined.)

\subsection{Comparison of \(\phi_S\) and \( \phi_S'\).}

Recall from \cite{ESL2} that \(\kh(L)\) and \(\kh'(L)\) are connected
by a spectral sequence. Indeed, the TQFT \(\A\) used to define
\(\kh\) is a graded TQFT, and this naturally gives rise to a grading
(known as the \(q\)-grading) on \(\ckh(D)\). Moreover, the TQFT \(\A'\) used
to define \(\ckh'(D)\) is a perturbation of \(D\): on the level of
groups, the two theories are isomorphic, and the maps induced by
cobordisms agree to lowest order in the \(q\)-grading. As a
consequence, the \(q\) grading defines a filtration on \(\ckh'(D)\),
which gives a spectral sequence \(E^i(D)\) converging to
\(\kh'(L)\). Moreover,  we have \(E^1(D) \cong \kh(L)\). This spectral
sequence behaves functorially with respect to cobordisms. To be
precise, we have 

\begin{lem}
\label{Lem:SSMap}
Let \( S \co L_1 \to L_2 \) be a cobordism with a fixed
decomposition into elementary cobordisms. Then for \(i \geq 1\), 
\( S \) induces a
morphism of spectral sequences  \( \phi_S^i \co E^i(L_1) \to
E^i(L_2) \) 
converging to \( \phi_{S}'\). If \( \overline{\phi}^1_{S}\co \kh(L_1)
\to \kh(L_2)\) is the induced map on filtered gradeds, then 
\( \overline{\phi}^1_{S} = \phi_S.\)
\end{lem}

\begin{proof}
This follows from standard properties of spectral sequences. Indeed, if
 \( f \co A \to B\) is a
map of filtered complexes, then there is an induced morphism of
spectral sequences \( f^i \co A^i \to B^i\) which converges to
\(f_* \). Moreover, if \(\overline{f}^0 \co \overline{A}^0 \to
 \overline{B}^0\) is the induced chain map on filtered gradeds, it is
 not difficult to see that \(\overline{f}^1 =
 \overline{f}^0_*\). 
 In the case when \(S\) is an elementary cobordism induced
by a Morse move, we are in precisely this situation with the map 
\(\phi \co \ckh(D_1) \to \ckh(D_2)\). The argument is very similar when
 \(S\) is an elementary cobordism associated to a Reidemeister move.
(See the proof of Theorem 1 in \cite{khg} for more details). Finally,
 the result for a general cobordism \(S\) follows by functoriality. 
\end{proof}

\section{Canonical generators}
\label{Sec:Gens}

Let \(L\) be a link represented by a planar diagram \(D\).
Given an orientation \(o\) of \(L\), Lee associates to it
a state \( \s_o \in \ckh'(D)\) and shows that \(\kh'(L)\) is freely
generated by \(\{[\s_o]\ts | \ts o \text{ is an orientation of } L\}\).
It turns out that there is a somewhat more natural way to normalize the
\(\s_o\)'s:

\begin{dfn}
Let \(w(o)\) be the writhe of the diagram \(D\)
endowed with the orientation \(o\), and let \(k(o)\) be the number of
circles in its oriented resolution. We define the rescaled canonical
generator associated to the orientation \(o\) by 
\begin{equation*}
\so = 2^{[w(o)-k(o)]/2} \s_o
\end{equation*}
\end{dfn}

To see why this choice of generators is a good one, we consider the
behavior of \( \so\) under the map induced by a Morse cobordism. 

\begin{prop}
\label{Prop:Map}
Let \(S:L_1 \to L_2\) be a cobordism with  no closed components. 
If \(o\) is an orientation on \(L_1\), then 
\begin{equation}
\label{Eq:Main}
\phi_S' (\sss_o) = 2^{-\chi(S)} \sum_I \pm \sss_{o_I}
\end{equation}
where the sum runs over all orientations on \(L_2\) compatible with
\(o\). 
\end{prop}

\begin{proof}

We  work our way up to the proof  in stages,
starting with the easiest possible case. 

\begin{lem}
Equation~\ref{Eq:Main} holds in the case where
 \(S: L_1 \to L_2\) is an elementary Morse cobordism.
\end{lem}

The proof is an easy calculation along the lines of the proof of
Proposition 4.1 in \cite{khg}. We leave its verification to the
reader. 

\begin{lem}
\label{Lem:ECob} 
Equation~\ref{Eq:Main} holds in the case where
 \(S:L_1 \to L_2\) is a composition of elementary Morse cobordisms
 with no closed component. 
\end{lem}

\begin{proof}
We induct on the number of elementary cobordisms in the
composition. The base case of one cobordism is covered by the lemma;
if there is more than one cobordism, we decompose \(S\) into the
composition of two cobordisms \(S_1: L_1 \to L_{1.5}\) and \(S_2 :
L_{1.5} \to L_2\) for which the statement is known to hold. Then we
compute 
\begin{equation*}
\phi_S' (\sss_o) = 2^{-\chi(S_1)-\chi(S_2)} \sum_{(o_{1.5},o_2)} \pm \sss_{o_2}
= 2^{-\chi(S)} \sum_{(o_{1.5},o_2)} \pm \sss_{o_2}
\end{equation*}
where the sum runs over pairs \((o_{1.5},o_2)\) such that \(o\) is
compatible with \(o_{1.5}\) and \(o_{1.5}\) is compatible with
\(o_2\). If this is the case, \(o_2\) is clearly compatible with
\(o_1\). Conversely, given \(o_2\) on \(L_2\) compatible with \(o\),
the fact that \(S\) has no closed components implies that there is a
unique orientation on \(S\) compatible with \(o\) and
\(o_2\). Thus \(o_2\) is associated with a unique compatible pair
\((o_{1.5},o_2)\). 
\end{proof}

To prove the proposition in general, we must check
that the generators \(\sss_o\) behave well under the isomorphisms
associated to Reidemeister moves.

\begin{lem}
\label{Lem:RMoves}
Suppose \(R_i \co D \to \tilde{D}\) is a Reidemeister move relating \(D\) to
another diagram \(\tilde{D}\) of \(L\), and let
\(\rho_{i*}' \co \kh'(D) \to \kh'(\tilde{D})\) be the corresponding
induced map. Let \(o\) be an orientation on \(D\) and \(\tilde{o}\) be
the correspoding orientation on \(\tilde{D}\). Then
\( \rho_{i*}'([\so]) = \pm [\overline{\s}_{\tilde{o}}]. \) 
\end{lem}

\begin{proof}
Essentially, this follows from the fact that the \( \rho_i\) are
defined in terms of maps induced by elementary cobordisms. 
Below, we give a more detailed argument for each Reidemeister move.

\vskip0.05in
\noindent{\it Reidemeister I Move:} Let \(\tilde{D}\) be a diagram
obtained from \(D\) by adding a left-hand curl, and let
\(\tilde{D}(*0)\) be the diagram obtained by giving this crossing the
 \(0\) resolution. Then the map
\(\rho_1 \co \ckh'(D_1) \to \ckh'(\tilde{D}(*0)) \subset
\ckh'(\tilde{D})\)  is given by \( \rho_i = \phi_{S_1}' -
\phi_{S_2}' \), where \(S_1\) is the obvious \(1\)-handle cobordism
from \(D\) to \(\tilde{D}(*0)\), and \(S_2\) is the product
cobordism connect summed with a trivial torus, followed by the
addition of a zero handle. Both \(S_1\) and \(S_2\) have \( \chi =
-1\). Using Lemma~\ref{Lem:ECob} (with a little direct computation
to get the signs right), we see that 
\begin{align*}
\rho_1(\sss_o)  & = 2^{1/2}[\sss_{o_1} - (\sss_{o_1}-\sss_{{o_2}})] \\
& = 2^{1/2} \sss_{{o_2}} \\
& =  \sss_{\tilde{o}}.
\end{align*}
Here \(o_1\) and \(o_2\) are the two orientations on \(\tilde{D}(*0)\)
compatible with \(o\) under \(S_2\). The last equality follows from
the fact that \( w(\tilde{D}) = w(\tilde{D}(*0))+1\). 

\begin{figure}
\includegraphics{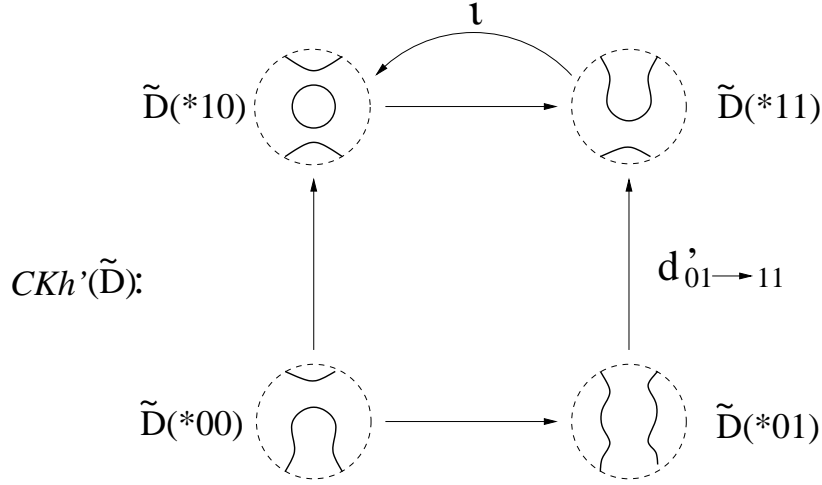}
\caption{\label{Fig:RTwo} Chain complex for the Reidemeister II move.}
\end{figure}

\vskip0.05in
\noindent{\it Reidemeister II Move:} Let \(\tilde{D}\) be obtained
from \(D\) by adding a pair of cancelling intersections. The complex
\( \ckh'(\tilde{D})\) is illustrated in Figure~\ref{Fig:RTwo}. The
part labeled \(\tilde{D}(*01)\) is naturally identified with \(
\ckh'(D)\), and for \(x\in \ckh'(D)\), we have \( \rho_2'(x) = x +
\alpha (x)\), where the map \( \alpha: \ckh'(\tilde{D}(*01)) \to
\ckh'(\tilde{D}(*10)) \) is  induced by a cobordism \(S \co 
\tilde{D}(*01) \to \tilde{D}(*10) \). \(S\) is 
the  composition of  two elementary
Morse cobordisms: from \(\tilde{D}(*01)\) to \( \tilde{D}(*11)\) by addition of
a 1-handle, and then from  \( \tilde{D}(*11)\) to \( \tilde{D}(*10)\) by
addition of a 0-handle. 

Fix an orientation \(o\) on \(D\). We consider two cases. First,
suppose that the two strands involved in the move have parallel
orientations, so that the oriented resolution of \(\tilde{D}\) is the
oriented resolution of
 \(\tilde{D}(*01)\). There is no orientation on \(\tilde{D}(*11)\)
compatible with \(o\), so \( \alpha (\sss_o) = 0 \). Since
\(w(o) = w(\tilde{o})\), it follows that \(\rho_2'(\sss_o) =
\sss_{\tilde{o}}\). 

Now suppose that the two strands involved in the move have opposite
orientations, so that the oriented resolution of \(\tilde{D}\) is the
oriented resolution of
 \(\tilde{D}(*10)\). In this case, there are two orientations
\(o_1\) and \(o_2\) on \(\tilde{D}(*01)\) compatible with
\(o\). Since \(\chi(S) = 0\), 
Corollary~\ref{Lem:ECob} implies that
 \( \alpha(\sss_o) = \pm \sss_{o_1} \pm \sss_{o_2} \). Thus we have
 \( \rho_2'(\sss_o) = \sss_o \pm \sss_{o_1} \pm
\sss_{o_2} \). Now one of  \(\sss_{o_1}\) or \(\sss_{o_2}\) is equal to
\(\sss_{\tilde{o}}\) --- without loss of generality, let us assume it
is \(o_1\). Then, as was observed in the proof of Proposition 2.3 of
\cite{khg}, the remaining term
\(\sss_o \pm \sss_{o_2} \) is exact in \( \ckh'(\tilde{D})\). Thus we
have \( \rho_2'([\sss_o]) = \pm [\sss_{\tilde{o}}]\), and the claim holds in
this case as well. 

\begin{figure}
\includegraphics{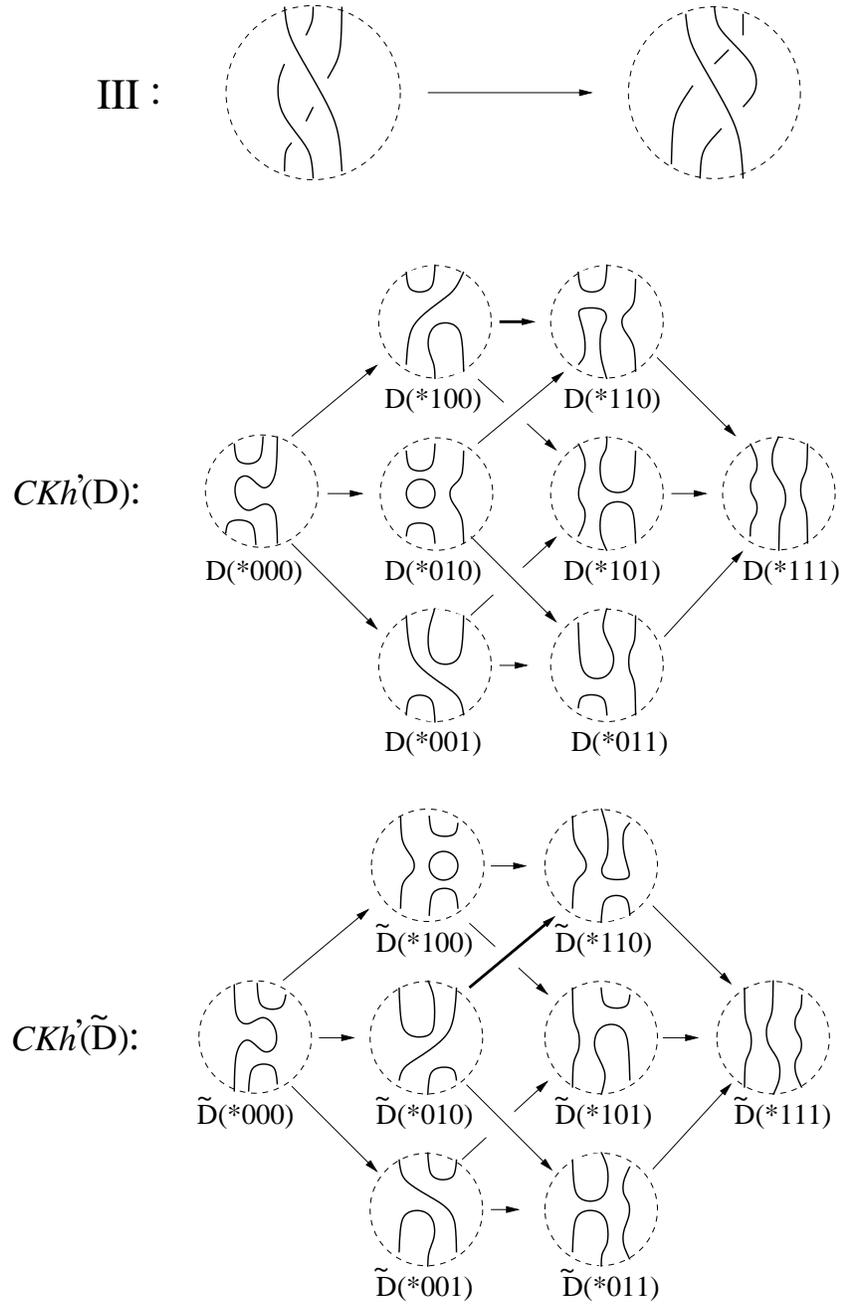}
\caption{\label{Fig:RThree} Chain complexes involved in
 the Reidemeister III Move}
\end{figure}
\vskip0.05in
\noindent{\it Reidemeister III Move:} Let \(D\) and \(\tilde{D}\) be
as shown in Figure~\ref{Fig:RThree}. In this case, we decompose the
complexes as shown in the figure. If \(\sss_o \) is contained in \(
D(*1)\) , the argument is easy. Indeed, there is a natural
identification \( i \co  D(*1) \to \tilde{D}(*1)\) , and we have
 \( \rho_3 (\sss_o) = i(\sss_o)) = \sss_{\tilde{o}} \). 

On the other
 hand, if \(\sss_o \) is contained in \(D(*0)\), the argument proceeds
 much as in the case of the Reidemeister II move. Given  \(\sss_o\),
 we first find a homologous element of the form \( x + \alpha (x)\) ,
 where \( x\in D(*100)\). Under \( \rho_3\), this element is mapped to 
\( x + \tilde{\alpha}(x)\) via the natural identification 
\( D(*100) = \tilde{D}(*010) \). As with the Reidemeister II
 move, there are two cases to consider. Either \( \sss_o \in D(*100)\),
 in which case \( \alpha (\sss_o) = \tilde{\alpha}(\sss_o) = 0 \), or 
\( \sss_o \in D(*010)\), in which case it is homologous to 
\( \pm \sss_{o_1}  \pm  \sss_o \pm \sss_{o_2} = \pm \sss_{o_1} +
 \alpha(\pm \sss_{o_1}) \). This in turn, maps to 
\( \pm \sss_{o_1} + \tilde{\alpha}(\pm \sss_{o_1})\), which is
 homologous to \( \pm \sss_{\tilde{o}}\).

\end{proof}

 Proposition~\ref{Prop:Map} now follows by combining
 Lemma~\ref{Lem:RMoves} with
 the proof of Lemma~\ref{Lem:ECob}. 
\end{proof}

We now specialize to the case where \(S \co K_1 \to K_2 \) is a
connected cobordism between two knots. In this case, there are two
orientations \(o\) and \(o'\) on \(K_1\). When there is no risk of
confusion, we denote the corresponding compatible
orientations on \(K_2\) by \(o\) and \(o'\) as well. Recall from
section 3 of \cite{khg} that the \(q\)-grading on \(\kh'(K_1)\) is
well-defined \(\text{mod} \ 4\), and that \(\sss_o+\sss_{o'}\) and
\(\sss_o-\sss_{o'}\) are homogenous elements whose \(q\)-gradings differ
by \(2 \ \text{mod} \ 4\). From the proposition, it follows that
 \(  \phi_S'(\sss_o-\sss_{o'})\) is a multiple of either
 \(\sss_o-\sss_{o'}\) or \(\sss_o+\sss_{o'}\); which one is determined by
 the \(\text{mod} \ 4\) \(q\)-grading of \(\phi_S'\). We have thus
 arrived at

\begin{cor}
Let \(S \co K_1 \to K_2 \) be a connected cobordism between two 
knots. Then if \(\chi(S) \equiv 0 \ \text{mod} \ 4\),
\begin{equation*}
 \phi_S' (\sss_o-\sss_{o'}) = \pm 2^{-\chi(S)} (\sss_o-\sss_{o'}) 
\end{equation*}
while if \(\chi(S) \equiv 2 \ \text{mod} \ 4\),
\begin{equation*}
 \phi_S' (\sss_o-\sss_{o'}) = \pm 2^{-\chi(S)} (\sss_o+\sss_{o'}).
\end{equation*}
\end{cor}

It is now easy to assemble the proof of the theorem stated in the
introduction. Indeed, suppose \( \Sigma \subset \R^4\) is a closed
surface of genus one. Then by choosing an appropriate Morse function
on \(\R^4\), we can decompose \(\Sigma\) into the union of a 0-handle,
a cobordism \( S \co U \to U\), and a 2-handle. It follows that 
\(n_\Sigma = |\epsilon(\phi_S(\x_+))|\), where \( \epsilon (\x_-) =
1\). 

We claim that in this case, we have 
\(\phi_S (\x_+)  = \phi_S' (\x_+) \). Indeed, by Lemma
\ref{Lem:SSMap}, we know that \(\phi_S\) is the map on filtered
gradeds associated to a morphism of spectral sequences which converges
to \(\phi_S'\). Since the knots in question are both the unknot, the
spectral sequences converge at the \(E^1\)  term. Thus we can write
\begin{equation*}
 \phi_S' (\x_+) = \phi_S (\x_+) + \psi(\x_+)
\end{equation*} 
where \(q(\psi(\x_+)) \geq q(\x_+) +\chi(S) + 4 = 3 \). Since \(\kh (U) \)
is trivial in \(q\)-grading \(3\) and higher, we must have \( \psi(\x_+) = 0
\). This proves the claim. 

\noindent We can now calculate
\begin{align*}
 \phi_S'(\x_+) 
& = \phi_S'(\frac{1}{2}(\s_o - {\s}_{o'})) \\
& = \phi_S'(2^{-1/2}(\so - \overline{\s}_{o'})) \\
& = \pm 2^{-\chi(S)/2} [2^{-1/2}(\so + \overline{\s}_{o'})] \\
& = \pm 2 [\frac{1}{2}(\s_o - {\s}_{o'})] \\
& = \pm 2 \x_-
\end{align*}
which  completes the proof. 
\qed

\bibliographystyle{plain}
\bibliography{../mybib}

\begin{thebibliography}{1}

\bibitem{DBN2}
D.~Bar-Natan.
\newblock {Khovanov's Homology for Tangles and Cobordisms}.
\newblock math.GT/0410495, 2004.

\bibitem{CarterSaito}
J.~S. Carter and M.~Saito.
\newblock Reidemeister moves for surface isotopies and their interpretation as
  moves to movies.
\newblock {\em J. Knot Theory Ramifications}, 2:251--284, 1993.

\bibitem{CSM}
J.~S. Carter, M.~Saito, and S.~Satoh.
\newblock {Ribbon-moves for 2-knots with 1-handles attached and
  Khovanov-Jacobsson numbers}.
\newblock math.GT/0407493, 2004.

\bibitem{Jacobsson}
M.~Jacobsson.
\newblock {An invariant of link cobordisms from Khovanov's homology theory}.
\newblock math.GT/0206303, 2002.

\bibitem{Khovanov}
M.~Khovanov.
\newblock A categorification of the {J}ones polynomial.
\newblock {\em Duke Math. J.}, 101:359--426, 2000.

\bibitem{Khovanov3}
M.~Khovanov.
\newblock {An invariant of tangle cobordisms}.
\newblock math.QA/0207264, 2002.

\bibitem{ESL2}
E.~S. Lee.
\newblock {K}hovanov's invariants for alternating links.
\newblock math.GT/0210213, 2002.

\bibitem{khg}
J.~Rasmussen.
\newblock Khovanov homology and the slice genus.
\newblock math.GT/04020131, 2004.

\bibitem{Tanaka}
K.~Tanaka.
\newblock {Khovanov-Jacobsson numbers and invariants of surface-knots derived
  from Bar-Natan's theory}.
\newblock math.GT/0502371, 2005.

\end{thebibliography}
\end{document}